\documentclass[12pt]{article}
\usepackage[centertags]{amsmath}
\usepackage{amsfonts}
\usepackage{amssymb}
\usepackage{latexsym}
\usepackage{amsthm}
\usepackage{newlfont}
\usepackage{graphicx}
\usepackage{listings}
\usepackage{booktabs}
\usepackage{abstract}
\date{}

\newlength{\defbaselineskip}
\newcommand{\setlinespacing}[1]%
           {\setlength{\baselineskip}{#1 \defbaselineskip}}

\newcommand{\disp}{\operatorname{disp}}
\newcommand{\N}{{\mathbb{N}}}

\newcommand{\actaqed}{\hfill $\actabox$}
{\medskip\noindent \textit{Proof of #1. }}%
{\actaqed \medskip}

\def\Di{{\mathcal D}}

\def\cB{\mathcal B}
\def\cC{{\mathcal C}}

\def \Tr{\mathcal T}

\def \cF{\mathcal F}
\def \cX{\mathcal X}

\def \cM{\mathcal M}
\def \cS{\mathcal S}
\def\R{{\mathbb R}}
\def\Z{\mathbb Z}

\def \T{\mathbb T}

\def \<{\langle}
\def\>{\rangle}

\def \La{\Lambda}
\def \ep{\epsilon}

\def \ff{\varphi}

\def\bt{\beta}

\def\la{\lambda}

\def \supp{\operatorname{supp}}
\def \sp{\operatorname{span}}

\def\ba{\mathbf a}
\def\bb{\mathbf b}
\def\bx{\mathbf x}
\def\by{\mathbf y}
\def\bz{\mathbf z}
\def\bk{\mathbf k}
\def\bu{\mathbf u}

\def\bm{\mathbf m}

\def\bp{\mathbf p}
\def\bs{\mathbf s}

\def\bN{\mathbf N}
\def\bW{\mathbf W}
\def\bH{\mathbf H}
\def\bB{\mathbf B}
\def\bK{\mathbf K}

\def\bt{\beta}

\newtheorem{Theorem}{Theorem}[section]
\newtheorem{Lemma}{Lemma}[section]
\newtheorem{Definition}{Definition}[section]
\newtheorem{Proposition}{Proposition}[section]
\newtheorem{Remark}{Remark}[section]

\newtheorem{Corollary}{Corollary}[section]
\newtheorem{Conjecture}{Conjecture}[section]
\numberwithin{equation}{section}

\newcommand{\be}{\begin{equation}}
\newcommand{\ee}{\end{equation}}

\begin{document}

\title{Connections between numerical integration, discrepancy, dispersion, and universal discretization}
\author{V.N. Temlyakov\thanks{University of South Carolina, Steklov Institute of Mathematics, and Lomonosov Moscow State University.  }}
\maketitle
\begin{abstract}
{The main goal of this paper is to provide a brief survey of recent results which connect together results from different areas of research. It is well known that numerical integration of functions with mixed 
smoothness is closely related to the discrepancy theory. We discuss this connection in detail and provide a general view of this connection. It was established recently that the new concept of {\it fixed volume discrepancy} is very useful in proving the upper bounds for the dispersion. Also, it was understood recently that point sets with small dispersion are very good for the universal discretization of the uniform norm of trigonometric polynomials. 
 }
\end{abstract}

\section{Introduction} 
\label{I} 

The problem of discretizing the $d$-dimensional unit cube $[0,1]^d$ is a fundamental problem of mathematics. Certainly, we should clarify what do we mean by {\it discretization}.
There are different ways of doing that. We can interpret $[0,1]^d$ as a compact set of $\R^d$ and use the idea of {\it covering numbers} (metric entropy). With such approach, for instance 
in the case of $\ell_\infty$ norm, we can find optimal coverings. For a given $n\in\N$ the regular grid with coordinates at the centers of intervals $[(k-1)/n,k/n]$, $k=1,\dots,n$, provide 
an optimal $\ell_\infty$ covering with the number of points $N=n^d$. Very often the unit cube 
$[0,1]^d$ plays the role of a domain, where smooth functions of $d$ variables are defined and we are interested in discretizing some continuous operations with these functions. A classical example of such a problem is the problem of numerical integration of functions. 
It turns out that the mentioned above regular grids are very far from being good economical discretizations of $[0,1]^d$ for numerical integration purposes. It is a fundamental problem of 
computational mathematics. Several areas of mathematical research are devoted to this problem: numerical integration, discrepancy, dispersion, sampling. Many nontrivial examples of good (in different sense) point sets are known (see, for instance, \cite{BC}, \cite{Mat}, \cite{NX}, \cite{RT}, \cite{TBook}, \cite{VT89}, \cite{VTbookMA}, \cite{DTU}). The main goal of this paper is to provide a brief survey of recent results which connect together results from different areas of research. It is well known that numerical integration of functions with mixed 
smoothness is closely related to the discrepancy theory. We discuss this connection in detail and provide a general view of this connection. It was established recently (see \cite{VT161}) that the new concept of {\it fixed volume discrepancy} is very useful in proving the upper bounds for the dispersion. Also, it was understood recently that point sets with small dispersion are very good for the universal discretization of the uniform norm of trigonometric polynomials (see \cite{VT160}). 

\section{Discrepancy as a special case of numerical integration}
\label{A}

We formulate the numerical integration problem in a general setting.  Numerical integration seeks good ways of approximating an integral
$$
\int_\Omega f(\bx)d\mu
$$
by an expression of the form
\be\label{A.1}
\La_m(f,\xi) :=\sum_{j=1}^m\la_jf(\xi^j),\quad \xi=(\xi^1,\dots,\xi^m),\quad \xi^j \in \Omega,\quad j=1,\dots,m. 
\ee
It is clear that we must assume that $f$ is integrable and defined at the points
 $\xi^1,\dots,\xi^m$. Expression (\ref{A.1}) is called a {\it cubature formula} $(\xi,\La)$ (if $\Omega \subset \R^d$, $d\ge 2$) or a {\it quadrature formula} $(\xi,\La)$ (if $\Omega \subset \R$) with knots $\xi =(\xi^1,\dots,\xi^m)$ and weights $\La:=(\la_1,\dots,\la_m)$. 
 
 Some classes of cubature formulas are of special interest. For instance, the Quasi-Monte Carlo cubature formulas, which have equal weights $1/m$, are important in applications. We use a special notation for these cubature formulas
 $$
 Q_m(f,\xi) :=\frac{1}{m}\sum_{j=1}^mf(\xi^j).
 $$
 
 The following class is a natural subclass of all cubature formulas. Let $B$ be a positive number and $Q(B,m)$ be the set of cubature formulas $\Lambda_m(\cdot,\xi)$ satisfying the additional condition
\be\label{A.2}
\sum_{\mu=1}^m |\lambda_\mu| \le B.
\ee
 
 For a function class $\bW$ we introduce a concept of error of the cubature formula $\La_m(\cdot,\xi)$ by
\be\label{A.3}
\La_m(\bW,\xi):= \sup_{f\in \bW} |\int_\Omega fd\mu -\La_m(f,\xi)|. 
\ee
The quantity $\La_m(\bW,\xi)$ is a classical characteristic of the quality of a given cubature formula $\La_m(\cdot,\xi)$. This setting is called {\it the worst case setting} in 
the Information Based Complexity. Typically, in approximation theory we study the behavior of 
the quantity $\La_m(\bW,\xi)$ for classes $\bW$ of smooth functions, in particular, for the unit balls of different spaces of smooth functions -- Sobolev, Nikol'skii, Besov spaces and spaces with mixed smoothness (see \cite{VTbookMA} and \cite{DTU}). The problem of finding optimal in the sense of order cubature formulas for a given class is of special importance. This means 
that we are looking for a cubature formula $\La_m^{opt}(\bW,\xi)$ such that
\be\label{A.4}
\La_m^{opt}(\bW,\xi) \asymp \inf_{\xi,\La}\La_m(\bW,\xi)=: \kappa_m(\bW).
\ee

We now describe some typical classes $\bW$, which are of interest in numerical integration and in discrepancy theory. We begin with a classical definition of discrepancy ("star discrepancy", $L_\infty$-discrepancy) of a point set $\xi := \{\xi^\mu\}_{\mu=1}^m\subset [0,1)^d$. 
Let $d\ge 2$ and $[0,1)^d$ be the $d$-dimensional unit cube. For convenience we sometimes use the notation $\Omega_d:=[0,1)^d$. For $\bx,\by \in [0,1)^d$ with $\bx=(x_1,\dots,x_d)$ and $\by=(y_1,\dots,y_d)$ we write $\bx < \by$ if this inequality holds coordinate-wise. For $\bx<\by$ we write $[\bx,\by)$ for the axis-parallel box $[x_1,y_1)\times\cdots\times[x_d,y_d)$ and define
$$
\cB:= \{[\bx,\by): \bx,\by\in [0,1)^d, \bx<\by\}.
$$

Introduce a class of special $d$-variate characteristic functions
$$
\chi^d := \{\chi_{[\mathbf 0,\bb)}(\bx):=\prod_{j=1}^d \chi_{[0,b_j)}(x_j),\quad b_j\in [0,1),\quad j=1,\dots,d\}
$$
where $\chi_{[a,b)}(x)$ is a univariate characteristic function of the interval $[a,b)$. 
The classical definition of discrepancy of a set $\xi$ of points $\{\xi^1,\dots,\xi^m\}\subset [0,1)^d$ is as follows
$$
D_\infty(\xi)  :=Q_m(\chi^d,\xi,\infty):= Q_m(\chi^d,\xi) = \max_{\bb\in [0,1)^d}\left|\prod_{j=1}^db_j -\frac{1}{m}\sum_{\mu=1}^m \chi_{[\mathbf 0,\bb)}(\xi^\mu)\right|.  
$$
The class $\chi^d$ is parametrized by the parameter $\bb\in [0,1)^d$. Therefore, we can define the $L_q$-discrepancy, $1\le q\le\infty$, of $\xi$ as follows
\be\label{A.5}
D_q(\xi) :=   Q_m(\chi^d,\xi,q):= \left\|\prod_{j=1}^db_j -\frac{1}{m}\sum_{\mu=1}^m \chi_{[\mathbf 0,\bb)}(\xi^\mu)\right\|_q,
\ee
where the $L_q$ norm is taken with respect to $\bb$ over the domain $\Omega_d$. 

We note that the fact that the class $\chi^d$ is parametrized by $\bb\in \Omega$ is a special 
important feature, which allows us to consider along with the worst case setting (\ref{A.4}) 
{\it the   average case setting} (\ref{A.5}).

\section{A brief history of results on classical discrepancy}
\label{B}

The first result on the lower bound for discrepancy was the following conjecture of van der Corput \cite{Co1} and \cite{Co2} formulated in 1935. Let $\xi^j\in[0,1]$, $j=1,2,\dots$, then we have
$$
\limsup_{m\to \infty}mD_\infty(\xi^1,\dots,\xi^m)  =\infty.
$$
This conjecture was proved by van Aardenne-Ehrenfest \cite{AE1} in 1945 (see also \cite{AE2}):
$$
\limsup_{m\to \infty}\frac{\log\log\log m}{\log\log m}mD_\infty(\xi^1,\dots,\xi^m) >0.
$$
We now list some classical lower estimates of discrepancy.
Let us denote
$$
D(m,d)_q :=\inf_{\xi} D_q(\xi),\quad \xi=\{\xi_j\}_{j=1}^m,\quad 1\le q\le \infty.
$$
In 1954 K. Roth \cite{Ro} proved that
\be\label{B.1}
D(m,d)_2 \ge C(d)m^{-1}(\log m)^{(d-1)/2}. 
\ee
In 1972 W. Schmidt \cite{Sch1} proved
\be\label{B.2}
D(m,2)_\infty \ge Cm^{-1}\log m . 
\ee
In 1977 W. Schmidt \cite{Sch} proved
\be\label{B.3}
D(m,d)_q \ge C(d,q)m^{-1}(\log m)^{(d-1)/2},\qquad 1<q\le \infty.  
\ee
In 1981 G. Hal{\' a}sz \cite{Ha} proved
\be\label{B.4}
D(m,d)_1 \ge C(d)m^{-1}(\log m)^{1/2}. 
\ee
The following conjecture has been formulated in \cite{BC} as an excruciatingly difficult great open problem.
\begin{Conjecture}\label{BC.1} We have for $d\ge 3$
$$
D(m,d)_\infty \ge C(d)m^{-1}(\log m)^{d-1}. 
$$
\end{Conjecture}
This problem is still open. Recently, D. Bilyk and M. Lacey \cite{BL} and D. Bilyk, M. Lacey, and A. Vagharshakyan \cite{BLV} proved
$$
D(m,d)_\infty \ge C(d)m^{-1}(\log m)^{(d-1)/2 + \delta(d)} 
$$
with some positive $\delta(d)$. 

For further historical discussion we refer the reader to surveys \cite{VT89}, \cite{Bi}, \cite{DTU}, and books \cite{BC}, \cite{Mat}, \cite{VTbookMA}. 

\section{Smooth discrepancy and numerical integration}
\label{C}

In the above definitions the function class $\chi^d$ with $d=1$ consists of characteristic functions, which have smoothness $1$ in the $L_1$ norm. In numerical integration it is natural to study 
function classes with arbitrary smoothness $r$. There are different generalizations of 
the above concept of discrepancy to the case of {\it smooth  discrepancy}. We discuss two of them here. In the definition of the first version of the $r$-discrepancy (see \cite{TBook}, \cite{VTbookMA}) instead of the characteristic function (this corresponds to $1$-discrepancy) we use the following function
\begin{align*}
B_r(\bx,\by)&:= \prod_{j=1}^d\bigl((r-1)!\bigr)^{-1}
(y_j - x_j )_+^{r-1},\\
  \bx,\by&\in\Omega_d,\qquad (a)_+ := \max (a,0).
\end{align*}
Denote
$$
\bB^{r,d}:=\{B_r(\bx,\by): \by \in \Omega_d\}.
$$
Then for a point set $\xi:=\{\xi^\mu\}_{\mu=1}^m$ of cardinality $m$ and weights $\Lambda:=\{\la_\mu\}_{\mu=1}^m$ we define the $r$-discrepancy of the pair $(\xi,\Lambda)$ by the formula
$$
D^r_q (\xi,\Lambda):=D(\xi,\Lambda,B_r,q):=\La_m(\bB^{r,d},\xi,q)
$$
\be\label{C.1}
 := \left \|\sum_{\mu=1}^{m}\lambda_{\mu}B_r (\xi^{\mu},\by)-
\prod_{j=1}^d (y_j^r /r!)\right\|_q  .
\ee

Consider the class
$ \dot{ \bW}_{p}^r  :=\bW^{B_r}_p$  consisting of the functions
$f(\bx)$ representable in the form
$$
f(\bx) =\int_{\Omega_d} B_r (\bx,\by)\varphi(\by) d \by,
\qquad \|\varphi\|_p \le 1.
$$
In connection with the definition of the class
$ \dot{ \bW}_{p}^r $ we remark here that
for the error of the cubature formula $(\xi,\Lambda)$ with
weights $\Lambda = (\lambda_1,\dots,\lambda_m)$ and knots
$\xi = (\xi^1,\dots,\xi^m)$ the following relation holds with  $p' := p/(p-1)$ 
\be\label{C.2}
\Lambda_m\bigl( \dot{ \bW}_{p}^r ,\xi\bigr)=\left \|\sum_{\mu=1}^{m}\lambda_{\mu}B_r (\xi^{\mu},\by)-
\prod_{j=1}^d (y_j^r /r!)\right\|_{p'} =
D^r_{p'} (\xi,\Lambda)  .
\ee
Thus, errors of numerical integration of classes $\dot{ \bW}_{p}^r$ are dual to the average errors of numerical integration of classes $\bB^{r,d}$. 

We now consider classes of periodic functions with mixed smoothness. For $\bx=(x_1,\dots,x_d)$ denote
$$
F_r(\bx) := \prod_{j=1}^d F_r(x_j),\quad F_r(x_j):= 1+2\sum_{k=1}^\infty k^{-r}\cos (kx_j-r\pi/2), 
$$
and
$$
\bW^r_p := \{f:f=\varphi\ast F_r,\quad \|\varphi\|_p \le 1\}.
$$
For $f\in \bW^r_p$ we denote $f^{(r)} :=\varphi$ where $\varphi$ is such that $f=\varphi\ast F_r$.
  In the case of integer $r$ the class $\bW^r_p$ is very close to the class of functions $f$, 
satisfying $\|f^{(r,\dots,r)}\|_p \le 1$, where $f^{(r,\dots,r)}$ is the mixed derivative of $f$ of order $rd$.  
\begin{Proposition}\label{CP.1} There exist two positive constants $C_1(d)$ and $C_2(d)$ such that for any $\La_m(\cdot,\xi)$ with a property $\sum_j \la_j =1$ we have
\be\label{C.4}
C_1(d)\La_m(\chi^d,\xi)\le \La_m(\bW^1_1,\xi) \le C_2(d)\La_m(\chi^d,\xi).  
\ee
\end{Proposition}
The reader can find the proof of Proposition \ref{CP.1} in \cite{VT89}. The following theorem is from \cite{VT89} (see also \cite{VTbookMA}, p. 250). 

\begin{Theorem}\label{CT.1} Let $1 \le p \le \infty$. Then for $r\in \N$
\be\label{C.5}
\kappa_m\bigl( \dot{ \bW}_p^r(\Omega_d)\bigr)\asymp
\kappa_m\bigl( \bW_p^r\bigr).
\ee
\end{Theorem}

We now proceed to the $r$-smooth discrepancy. The classical definition of discrepancy of a set $\xi$ of points $\{\xi^1,\dots,\xi^m\}\subset [0,1)^d$   is equivalent within multiplicative constants, which may only depend on $d$, to the following definition
\be\label{C.6}
D^1(\xi):=  \sup_{B\in\cB}\left|vol(B)-\frac{1}{m}\sum_{\mu=1}^m \chi_B(\xi^\mu)\right|,
\ee
where for $B=[\ba,\bb)\in \cB$ we denote $\chi_B(\bx):= \prod_{j=1}^d \chi_{[a_j,b_j)}(x_j)$.   Moreover, we consider the following optimized version of $D^1(\xi)$
\be\label{C.7}
D^{1,o}(\xi):= \inf_{\la_1,\dots,\la_m} \sup_{B\in\cB}\left|vol(B)- \sum_{\mu=1}^m \la_\mu \chi_B(\xi^\mu)\right|.
\ee

We say that a univariate function $f$ has smoothness $1$ in $L_1$ if $\|\Delta_t f\|_1\le C|t|$, where $\Delta_tf(x):=f(x)-f(x+t)$. In case $\|\Delta^r_t f\|_1 \le C|t|^r$, where $\Delta^r_t:= (\Delta_t)^r$, we say that $f$ has smoothness $r$ in $L_1$.
In the definition of $D^1(\xi)$ and $D^{1,o}(\xi)$ -- the $1$-smooth discrepancy -- we use as a building block the univariate characteristic function.   In numerical integration $L_1$-smoothness of a function plays an important role. A characteristic function of an interval has smoothness $1$ in the $L_1$ norm. This is why we call the corresponding discrepancy characteristics the $1$-smooth discrepancy. In the definition of $D^2(\xi)$,
$D^{2,o}(\xi)$, $D^2(\xi,V)$, and $D^{2,o}(\xi,V)$ (see below and \cite{VT161}) we use the hat function 
$h_{[-u,u)}(x) =u-|x|$ for $|x|\le u$ and $h_{[-u,u)}(x) =0$ for $|x|\ge u$ instead of the characteristic function $\chi_{[-u/2,u/2)}(x)$. Function $h_{[-u,u)}(x)$ has smoothness $2$ in $L_1$. This fact gives the corresponding name. Note that
$$
h_{[-u,u)}(x) = \chi_{[-u/2,u/2)}(x) \ast \chi_{[-u/2,u/2)}(x),
$$
where
$$
f(x)\ast g(x) := \int_\R f(x-y)g(y)dy.
$$
Now, for $r=1,2,3,\dots$ we inductively define
$$
h^1(x,u):= \chi_{[-u/2,u/2)}(x),\qquad h^2(x,u):= h_{[-u,u)}(x),
$$
$$
h^r(x,u) := h^{r-1}(x,u)\ast h^1(x,u),\qquad r=3,4,\dots.
$$
Then $h^r(x,u)$ has smoothness $r$ in $L_1$ and has support $(-ru/2,ru/2)$. 
Represent a box $B\in\cB$ in the form
$$
B= \prod_{j=1}^d [x^0_j-ru_j/2,x^0_j+ru/2)
$$
and define
$$
h^r_B(\bx):= h^r(\bx,\bx^0,\bu):=\prod_{j=1}^d h^r(x_j-x^0_j,u_j).
$$

In \cite{VT161} we modified definitions (\ref{C.6}) and (\ref{C.7}), replacing the characteristic function $\chi_B$ by a smoother hat function $h^r_B$.  

The $r$-smooth discrepancy is now defined as
\be\label{C.8}
D^r(\xi):=  \sup_{B\in\cB}\left|\int h^r_B(\bx)d\bx-\frac{1}{m}\sum_{\mu=1}^m h^r_B(\xi^\mu)\right|
\ee
and its optimized version as 
\be\label{C.9}
D^{r,o}(\xi):=  \inf_{\la_1,\dots,\la_m}\sup_{B\in\cB}\left|\int h^r_B(\bx)d\bx- \sum_{\mu=1}^m \la_\mu h^r_B(\xi^\mu)\right|.
\ee
Note that the known concept of $r$-discrepancy   (see, for instance, \cite{TBook}, \cite{VT89}, and above in this section) is close to the concept of $r$-smooth discrepancy.

It is more convenient for us to consider the average setting in the periodic case. For a function $f\in L_1(\R^d)$ with a compact support we define its periodization $\tilde f$ as follows
$$
\tilde f(\bx) := \sum_{\bm\in \Z^d} f(\bm+\bx).
$$
For each $\bz\in [0,1)^d$ and $\bu \in (0,\frac{1}{2}]^d$ consider a periodization of function $h^r(\bx,\bz,\bu)$ 
in $\bx$ with period $1$ in each variable  
$\tilde h^r(\bx,\bz,\bu)$. Consider the class of periodic $r$-smooth hat functions
$$
\bH^{r,d}:=\{\tilde h^r(\bx,\bz,\bu): \bz\in[0,1)^d;\bu\in(0,1/2]^d\}.
$$

Define the corresponding {\it periodic $r$-smooth discrepancy}   as follows
$$
\tilde D^{r}_\infty(\xi,\La):=  \La_m(\bH^{r,d},\xi)
$$
\be\label{2.4}
   = \sup_{\bz\in[0,1)^d;\bu\in(0,1/2]^d}\left|\int_{[0,1)^d} \tilde h^r(\bx,\bz,\bu)d\bx- \sum_{\mu=1}^m \la_\mu \tilde h^r(\xi^\mu,\bz,\bu)\right|.
\ee
For $1\le p_1,p_2\le \infty$,  define the corresponding {\it periodic $r$-smooth $L_\bp$-discrepancy}, which also can be called {\it Weyl $r$-smooth $L_\bp$-discrepancy} (\cite{W}, \cite{Lev7}), as follows (see \cite{VT163} for the case $p=\infty$)
\be\label{2.5}
  \tilde D^{r}_{p_1,p_2}(\xi,\La):=
  \left\|\|\int_{[0,1)^d} \tilde h^r(\bx,\bz,\bu)d\bx- \sum_{\mu=1}^m \la_\mu \tilde h^r(\xi^\mu,\bz,\bu)\|_{p_1}\right\|_{p_2}
\ee
where the $L_{p_1}$ norm is taken with respect to $\bz$ over the unit cube $[0,1)^d$ and the $L_{p_2}$ norm is taken with respect to $\bu$ over the cube $(0,1/2]^d$. In the definition of $\tilde D^{r}_{p_1,p_2}(\xi,\La)$ parameters $\bz$ and $\bu$ play different roles. The most important parameter is $\bu$ -- it controls the shape of supports of the corresponding hat functions. It seems like the most natural value for parameter $p_2$ is $\infty$. In this case we obtain bounds uniform with respect to the shape and the size of supports of hat functions. 

\section{Lower estimates for the smooth discrepancy}
\label{D}

We now present the results on the lower estimates for the $r$-discrepancy. As above for a point set $\xi:=\{\xi^\mu\}_{\mu=1}^m$ of cardinality $m$ and weights $\Lambda:=\{\la_\mu\}_{\mu=1}^m$ we define the $r$-discrepancy of the pair $(\xi,\Lambda)$ by the formula
$$
D^r_q (\xi,\Lambda):=D(\xi,\Lambda,B_r,q):=\La_m(\bB^{r,d},\xi,q)
$$
\be\label{D.1}
 := \left \|\sum_{\mu=1}^{m}\lambda_{\mu}B_r (\xi^{\mu},\by)-
\prod_{j=1}^d (y_j^r /r!)\right\|_q  .
\ee
We
 denote
$$
D^r_q(m,d) := \inf_{\xi}D^r_q(\xi,(1/m,\dots,1/m))
$$
where $D^r_q(\xi,\Lambda)$ is defined in (\ref{D.1}) and also denote
$$
 D^{r,o}_q(m,d) := \inf_{\xi,\Lambda}D^r_q(\xi,\Lambda).
$$
It is clear that 
$$
 D^{r,o}_q(m,d)\le D^r_q(m,d).
$$
The first result on estimating the $r$-discrepancy was obtained in 1985 by V.A. Bykovskii \cite{By}
\be\label{D.2}
D_2^{r,o}(m,d) \ge C(r,d)m^{-r}(\log m)^{(d-1)/2}. 
\ee
This result is a generalization of Roth's result (\ref{B.1}).
The generalization of Schmidt's result (\ref{B.3}) was obtained by the author in 1990 (see \cite{VT43})  
\be\label{D.3}
D^{r,o}_q(m,d) \ge C(r,d,q)m^{-r}(\log m)^{(d-1)/2}, \qquad 1<q\le \infty. 
\ee
In 1994 (see \cite{VT50}) the author proved the lower bounds in the case of weights $\Lambda$ satisfying an extra condition (\ref{A.2}).
\begin{Theorem}\label{DT.1} Let $B$ be a positive number. For any points $\xi^1,\dots,\xi^m \subset \Omega_d$ and any weights $\Lambda =(\lambda_1,\dots,\lambda_m)$ satisfying the condition
\be\label{D.4}
\sum_{\mu=1}^m|\lambda_\mu| \le B
\ee
we have for even integers $r$
$$
D^r_\infty(\xi,\Lambda) \ge C(d,B,r)m^{-r}(\log m)^{d-1}
$$
with a positive constant $C(d,B,r)$.
\end{Theorem}
This result encouraged us to formulate the following generalization of the Conjecture \ref{BC.1} (see \cite{VT89}).
\begin{Conjecture}\label{DC.1} For all $d,r\in \N$ we have
$$
D^{r,o}_\infty(m,d) \ge C(r,d)m^{-r}(\log m)^{d-1}. 
$$
\end{Conjecture}

We now proceed to the $r$-smooth $L_\bp$-discrepancy. The first lower bound for such discrepancy was obtained in the case $\bp=\infty$ under an extra condition (\ref{D.4}) on the weights (see \cite{VT163}). Here is the corresponding result from \cite{VT163}.
\begin{Theorem}\label{DT.2}   For any points $\xi^1,\dots,\xi^m \subset \Omega_d$ and  weights $\Lambda =(\lambda_1,\dots,\lambda_m)$ satisfying condition (\ref{D.4})
we have for even integers $r$
$$
\tilde D^r_\infty(\xi,\Lambda) \ge C(d,B,r)m^{-r}(\log m)^{d-1}
$$
with a positive constant $C(d,B,r)$.
\end{Theorem}
Denote as above
$$
 \tilde D^{r,o}_\bp(m,d) := \inf_{\xi,\Lambda}\tilde D^r_\bp(\xi,\Lambda).
$$
Theorem \ref{DT.2} supports the following conjecture.
\begin{Conjecture}\label{DC.2} For all $d,r\in \N$ we have
$$
\tilde D^{r,o}_\infty(m,d) \ge C(r,d)m^{-r}(\log m)^{d-1}. 
$$
\end{Conjecture}

The following theorem is from \cite{VT164}.

 \begin{Theorem}\label{DT.3} Let $r\in\N$. Then for any $(\xi,\La)$ we have
$$
    \tilde D^{r}_{2,2}(\xi,\La) \geq   C(r,d)   m^{-r}(\log   m)^{(d-1)/2},   \qquad
C(r,d)>0.
$$
\end{Theorem}

Theorem \ref{DT.3} gives the following lower bound for $r\in\N$ and $\bp\ge \mathbf 2$
\be\label{D.5}
 \tilde D^{r,o}_\bp(m,d) \ge C(r,d)m^{-r}(\log m)^{(d-1)/2}. 
\ee
The lower bound (\ref{D.5}) is different from the lower bound from Theorem \ref{DT.2}. 
However, the following Proposition \ref{DP.1} (see \cite{VT164}) shows that this bound is sharp in case $\bp=\mathbf 2$.

\begin{Proposition}\label{DP.1} For $r\in \N$ there exists a cubature formula $(\xi,\La)$
such that 
$$
    \tilde D^{r}_{2,\infty}(\xi,\La) \le   C(r,d)   m^{-r}(\log   m)^{(d-1)/2},   \qquad
C(r,d)>0.
$$
\end{Proposition}
 
Under stronger assumption on $r$, namely, assuming that $r$ is an even number, we obtained in \cite{VT164} a stronger than (\ref{D.5}) lower bound. 

\begin{Theorem}\label{DT.4} Let $r\in\N$ be an even number. Then for any $(\xi,\La)$ we have for $1<p<\infty$
$$
    \tilde D^{r}_{p,1}(\xi,\La) \geq   C(r,d,p)   m^{-r}(\log   m)^{(d-1)/2},   \qquad
C(r,d,p)>0.
$$
\end{Theorem}

Theorem \ref{DT.4} gives that for even $r$
 for any   $1<p<\infty$
$$
    \tilde D^{r,o}_{p,1}(m,d) \geq   C(r,d,p)   m^{-r}(\log   m)^{(d-1)/2},   \qquad
C(r,d,p)>0.
$$
The following 
 result from \cite{VT164} is an extension of Proposition \ref{DP.1}. 
 \begin{Proposition}\label{DP.2} For $r\in \N$ and $1<p<\infty$ there exists a cubature formula $(\xi,\La)$
such that 
$$
    \tilde D^{r}_{p,\infty}(\xi,\La) \le   C(r,p,d)   m^{-r}(\log   m)^{(d-1)/2}.
$$
\end{Proposition}

Proposition \ref{DP.2} shows that the above lower bound is sharp. Moreover, it shows that for $r$ even we have for all $1<p_1<\infty$ and $1\le p_2 \le \infty$
\be\label{D.6}
\tilde D^{r,o}_{p_1,p_2}(m,d) \asymp       m^{-r}(\log   m)^{(d-1)/2} .
\ee

\section{Fixed volume discrepancy}
\label{E}

Along with $D^r(\xi)$ and $D^{r,o}(\xi)$ we consider a more refined quantity -- {\it $r$-smooth fixed volume discrepancy} -- defined as follows
\be\label{E.1}
D^r(\xi,V):=  \sup_{B\in\cB:vol(B)=V}\left|\int h_B^r(\bx)d\bx-\frac{1}{m}\sum_{\mu=1}^m h_B^r(\xi^\mu)\right|;
\ee
\be\label{E.2}
D^{r,o}(\xi,V):=  \inf_{\la_1,\dots,\la_m}\sup_{B\in\cB:vol(B)=V}\left|\int h_B^r(\bx)d\bx- \sum_{\mu=1}^m \la_\mu h_B^r(\xi^\mu)\right|.
\ee
Clearly,
$$
D^r(\xi) = \sup_{V\in(0,1]} D^r(\xi,V).
$$

We begin with the case $d=2$. It is well known that the Fibonacci cubature formulas are optimal in the sense of order for numerical integration of different kind of smoothness classes of functions of two variables (see \cite{TBook}, \cite{VTbookMA}, \cite{DTU}). We present a result from \cite{VT161}, which shows that the Fibonacci point set has good fixed volume discrepancy.

Let $\{b_n\}_{n=0}^{\infty}$, $b_0=b_1 =1$, $b_n = b_{n-1}+b_{n-2}$,
$n\ge 2$, --
be the Fibonacci numbers. Denote the $n$th {\it Fibonacci point set} by
$$
\cF_n:= \left\{(\mu/b_n,\{\mu b_{n-1} /b_n \}),\, \mu=1,\dots,b_n\right\}.
$$
In this definition $\{a\}$ is the fractional part of the number $a$.   The cardinality of the set $\cF_n$ is equal to $b_n$. In \cite{VT161} we proved 
the following upper bound.

\begin{Theorem}\label{ET.1} Let $d=2$, $r\ge2$. There exists a   constant $c(r)>0$ such that for any $V\ge V_0:= c(r)/b_n$ we have for all $B\in\cB$, $vol(B)=V$
\be\label{E.3}
\left|b_n^{-1}\sum_{\mu=1}^{b_n} h^r_B(\mu/b_n,\{\mu b_{n-1}/b_n\}) - \hat h^r_B(\mathbf 0)\right| \le C(r)\log(2V/V_0)/b_n^r.
\ee
\end{Theorem}

  Theorem \ref{ET.1} 
provides the following inequalities for the Fibonacci point sets $\cF_n$ in case $r\ge 2$
$$
D^{r,o}(\cF_n,V) \le D^{r}(\cF_n,V) \le C(r)(\log(2V/V_0))/b_n^r,\qquad V\ge V_0.
$$ 

We now proceed to the case $d\ge 3$. It is well known that the {\it Frolov point sets} are very good for numerical integration of  smoothness classes of functions of several variables (see \cite{Fro1}, \cite{TBook}, \cite{VT89}, \cite{VTbookMA}, \cite{DTU}, \cite{MUll}). Theorem \ref{ET.2} below, which was proved in \cite{VT161}, shows that the Frolov point sets have good fixed volume discrepancy. Construction of the Frolov point sets 
is more involved than the construction of the Fibonacci point sets. We begin with a description of the Frolov point sets. The following lemma plays a fundamental role in the
construction of such point sets (see \cite{TBook}, Ch.4, \S4 or \cite{VTbookMA}, Ch.6, S.6.7 for its proof).

\begin{Lemma}\label{EL.1} There exists a matrix $A$ such that the lattice
$L(\mathbf m) = A\mathbf m$
$$
 L(\mathbf m) =
\begin{bmatrix}
L_1(\mathbf m)\\
\vdots\\
L_d(\mathbf m)
\end{bmatrix},
$$
where $\mathbf m$ is a (column)
vector with integer coordinates,
has the following properties

{$1^0$}. $\qquad \left |\prod_{j=1}^d L_j(\mathbf m)\right|\ge 1$
for all $\mathbf m \ne \mathbf 0$;

{$2^0$} each parallelepiped $P$ with volume $|P|$
whose edges are parallel
to the coordinate axes contains no more than $|P| + 1$ lattice
points.
\end{Lemma}

Let $a > 1$ and $A$ be the matrix from Lemma \ref{EL.1}. We consider the
cubature formula
$$
\Phi(a,A)(f) := \bigl(a^d |\det A|\bigr)^{-1}\sum_{\mathbf m\in\Z^d}f
\left (\frac{(A^{-1})^T\mathbf m}{a}\right)
$$
for $f$ with compact support.   

We call the {\it Frolov point set} the following set associated with the matrix $A$ and parameter $a$
$$
\cF(a,A) := \left\{\left (\frac{(A^{-1})^T\mathbf m}{a}\right)\right\}_{\bm\in\Z^d}\cap [0,1]^d =: \{z^\mu\}_{\mu=1}^N.
$$
 Clearly, the number $N=|\cF(a,A)|$ of points of this
set does not exceed $C(A)a^d$.

\begin{Theorem}\label{ET.2} Let $r\ge 2$. There exists a constant $c(d,A,r)>0$ such that for any $V\ge V_0:= c(d,A,r)a^{-d}$ we have for all $B\in\cB$, $vol(B)=V$,
\be\label{E.4}
|\Phi(a,A)(h^r_B) - \hat h^r_B(\mathbf 0)| \le C(d,A,r)a^{-rd} (\log(2V/V_0))^{d-1}.
\ee
\end{Theorem}

\begin{Corollary}\label{EC.1} For $r\ge2$ there exists a constant $c(d,A,r)>0$ such that for any $V\ge V_0:= c(d,A,r)a^{-d}$ we have
\be\label{E.5}
D^{r,o}(\cF(a,A),V) \le C(d,A,r)a^{-rd} (\log(2V/V_0))^{d-1}.
\ee
\end{Corollary}

The following technical Lemma \ref{EL.2} played the main role in the proofs of Theorems \ref{ET.1} and \ref{ET.2}. Lemma \ref{EL.2} might be of interest by itself. Consider
$$
\sigma^r(v,\bu):= \sum_{\|\bs\|_1=v}\prod_{j=1}^d \min\left((2^{s_j}u_j)^{r/2},\frac{1}{(2^{s_j}u_j)^{r/2}}\right),\quad v\in\N_0.
$$
Denote
$$
pr(\bu,d) := \prod_{j=1}^d u_j.
$$
 
 \begin{Lemma}\label{EL.2} Let $v\in \N_0$ and $\bu\in \R^d_+$. Then we have 
the following inequalities.

(I) Under condition $2^vpr(\bu,d)\ge 1$ we have
\be\label{E.6}
\sigma^r(v,\bu) \le C(d) \frac{\left(\log(2^{v+1}pr(\bu,d))\right)^{d-1}}{(2^vpr(\bu,d))^{r/2}}.
\ee

(II) Under condition $2^vpr(\bu,d)\le 1$ we have
\be\label{E.7}
\sigma^r(v,\bu) \le C(d) (2^vpr(\bu,d))^{r/2} \left(\log\frac{2}{2^{v}pr(\bu,d)}\right)^{d-1}.
\ee

\end{Lemma}

In \cite{VT163} we extended Theorem \ref{ET.2} and Corollary \ref{EC.1} to the periodic case. For that we need to modify the set $\cF(a,A)$ and the cubature formula 
$\Phi(a,A)$. For $\by\in \mathbb R^d$ denote $\{\by\} := (\{y_1\},\dots,\{y_d\})$, where for $y\in \mathbb R$ notation $\{y\}$ means the fractional part of $y$. For given $a$ and $A$ denote
$$
\eta:= \{\eta^\mu\}_{\mu=1}^m := \left\{\left (\frac{(A^{-1})^T\mathbf m}{a}\right)\right\}_{\bm\in\Z^d}\cap [-1/2,3/2)^d
$$
and
\be\label{E.8}
\xi:=\{\xi^\mu\}_{\mu=1}^m := \left\{\{\eta^\mu\}\right\}_{\mu=1}^m.
\ee
Clearly, $m\le C(A)a^d$. Next, let $w(t)$ be infinitely differentiable on $\mathbb R$ function 
with the following properties
\be\label{E.9}
\supp(w) \subset (-1/2,3/2)\quad \text{and}\quad \sum_{k\in\Z} w(t+k) =1.
\ee
Denote $w(\bx):= \prod_{j=1}^d w(x_j)$. Then for $f(\bx)$ defined on $[0,1)^d$ we 
consider the cubature formula
$$
\Phi(a,A,w)(f) := \sum_{\mu=1}^m w_\mu f(\xi^\mu),\qquad w_\mu := w(\eta^\mu) .
$$
 In \cite{VT163} we proved the following analogs of Theorem \ref{ET.2} and Corollary \ref{EC.1}.
 
 \begin{Theorem}\label{ET.3} Let $r\ge 2$. There exists a constant $c(d,A,r)>0$ such that for any $v\ge v_0:= c(d,A,r)a^{-d}$ we have for all $\bu\in (0,1/2]^d$, $pr(\bu,d)=v$, and $\bz\in[0,1)^d$
$$
|\Phi(a,A,w)(\tilde h^r(\cdot,\bz,\bu)) - \hat{\tilde h}^r(\mathbf 0,\bz,\bu)| \le C(d,A,r,w)a^{-rd} (\log(2v/v_0))^{d-1}.
$$
\end{Theorem}

\begin{Corollary}\label{EC.2} For $r\ge2$ there exists a constant $c(d,A,r)>0$ such that for any $v\ge v_0:= c(d,A,r)a^{-d}$ we have for the point set $\xi$ defined by (\ref{E.8})
$$
\tilde D^{r,o}(\xi,v) \le C(d,A,r)a^{-rd} (\log(2v/v_0))^{d-1}.
$$
\end{Corollary}

In particular, Theorem \ref{ET.3} implies that the $r$-smooth periodic discrepancy
$$
\tilde D^{r,o}_m :=
$$
\be\label{E.10}
   \inf_{\la_1,\dots,\la_m}\sup_{\bz\in[0,1)^d;\bu\in(0,1/2]^d}\left|\int_{[0,1)^d} \tilde h^r(\bx,\bz,\bu)d\bx- \sum_{\mu=1}^m \la_\mu \tilde h^r(\xi^\mu,\bz,\bu)\right|
\ee
satisfies the bound (for $r\in\N$, $r\ge 2$)
\be\label{E.11}
\tilde D^{r,o}_m   \le C(d,r)m^{-r} (\log m)^{d-1}.
\ee
  Theorem \ref{DT.2} shows that the bound (\ref{E.11}) cannot be improved for a natural 
class of weights $\la_1,\dots,\la_m$ used in the optimization procedure in the definition 
of $\tilde D^{r,o}_m$, namely, for weights, satisfying
$$
\sum_{\mu=1}^m |\la_\mu| \le B.
$$

 \section{Dispersion}
 \label{F}
 
 We remind the definition of 
dispersion. Let $d\ge 2$ and $[0,1)^d$ be the $d$-dimensional unit cube. As above for $\bx,\by \in [0,1)^d$ with $\bx=(x_1,\dots,x_d)$ and $\by=(y_1,\dots,y_d)$ we write $\bx < \by$ if this inequality holds coordinate-wise. For $\bx<\by$ we write $[\bx,\by)$ for the axis-parallel box $[x_1,y_1)\times\cdots\times[x_d,y_d)$ and define
$$
\cB:= \{[\bx,\by): \bx,\by\in [0,1)^d, \bx<\by\}.
$$
For $n\ge 1$ let $T$ be a set of points in $[0,1)^d$ of cardinality $|T|=n$. The volume of the largest empty (from points of $T$) axis-parallel box, which can be inscribed in $[0,1)^d$, is called the dispersion of $T$:
$$
\text{disp}(T) := \sup_{B\in\cB: B\cap T =\emptyset} vol(B).
$$
An interesting extremal problem is to find (estimate) the minimal dispersion of point sets of fixed cardinality:
$$
\text{disp*}(n,d) := \inf_{T\subset [0,1)^d, |T|=n} \text{disp}(T).
$$
It is known that 
\be\label{F.1}
\text{disp*}(n,d) \le C^*(d)/n.
\ee
Inequality (\ref{F.1}) with $C^*(d)=2^{d-1}\prod_{i=1}^{d-1}p_i$, where $p_i$ denotes the $i$th prime number, was proved in \cite{DJ} (see also \cite{RT}). The authors of \cite{DJ} used the Halton-Hammersly set of $n$ points (see \cite{Mat}). Inequality (\ref{F.1}) with $C^*(d)=2^{7d+1}$ was proved in 
\cite{AHR}. The authors of \cite{AHR}, following G. Larcher, used the $(t,r,d)$-nets (see \cite{NX}, \cite{Mat} for results on $(t,r,d)$-nets and Definition \ref{GD.1} below for the definition). 

It was demonstrated in \cite{VT161} how good upper bounds on fixed volume discrepancy can be used for proving good upper bounds for dispersion. This fact was one of the motivation for studying the fixed volume discrepancy. 
Theorem \ref{FT.1} below was derived from Theorem \ref{ET.1} (see \cite{VT161}). The upper bound in Theorem \ref{FT.1} combined with 
the trivial lower bound shows that the Fibonacci point set provides optimal rate of decay for the dispersion. 

\begin{Theorem}\label{FT.1} There is an absolute constant $C$ such that for all $n$ we have
\be\label{F.2}
\disp (\cF_n) \le C/b_n.
\ee
\end{Theorem}

The following Theorem \ref{FT.2} was derived in \cite{VT161} from Theorem \ref{ET.2}.

\begin{Theorem}\label{FT.2} Let $A$ be a matrix from Lemma \ref{EL.1}. There is a  constant $C(d,A)$, which may only depend on $A$ and $d$, such that for all $a$ we have
\be\label{F.3}
\disp(\cF(a,A)) \le C(A,d)a^{-d}.
\ee
\end{Theorem}

The reader can find further recent results on dispersion in \cite{Rud}, \cite{Sos}, and \cite{Ull}. 

\section{Universal discretization of the uniform norm}
\label{G}

In this section we demonstrate an application of results on dispersion from Section \ref{F}  to the problem of universal discretization. For a more detailed discussion of 
universality in approximation and learning theory we refer the reader to \cite{Tem16}, 
\cite{TBook}, \cite{VT89}, \cite{VTbookMA}, \cite{DTU}, \cite{VT160}, \cite{GKKW},
\cite{BCDDT}, \cite{VT113}. We remind the discretization problem setting, which we plan to discuss (see \cite{VT158} and \cite{VT159}). 

{\bf Marcinkiewicz problem.} Let $\Omega$ be a compact subset of $\R^d$ with the probability measure $\mu$. We say that a linear subspace $X_N$ (usually $N$ stands for the dimension of $X_N$) of the $L_q(\Omega)$, $1\le q < \infty$, admits the Marcinkiewicz-type discretization theorem with parameters $m$ and $q$ if there exist a set $\{\xi^\nu \in \Omega, \nu=1,\dots,m\}$ and two positive constants $C_j(d,q)$, $j=1,2$, such that for any $f\in X_N$ we have
\be\label{G.1}
C_1(d,q)\|f\|_q^q \le \frac{1}{m} \sum_{\nu=1}^m |f(\xi^\nu)|^q \le C_2(d,q)\|f\|_q^q.
\ee
In the case $q=\infty$ we define $L_\infty$ as the space of continuous on $\Omega$ functions and ask for 
\be\label{G.2}
C_1(d)\|f\|_\infty \le \max_{1\le\nu\le m} |f(\xi^\nu)| \le  \|f\|_\infty.
\ee
We will also use a brief way to express the above property: the $\cM(m,q)$ theorem holds for  a subspace $X_N$ or $X_N \in \cM(m,q)$. 

{\bf Universal discretization problem.} This problem is about finding (proving existence) of 
a set of points, which is good in the sense of the above Marcinkiewicz-type discretization 
for a collection of linear subspaces (see \cite{VT160}). We formulate it in an explicit form. Let $\cX_N:= \{X_N^j\}_{j=1}^k$ be a collection of linear subspaces $X_N^j$ of the $L_q(\Omega)$, $1\le q \le \infty$. We say that a set $\{\xi^\nu \in \Omega, \nu=1,\dots,m\}$ provides {\it universal discretization} for the collection $\cX_N$ if, in the case $1\le q<\infty$, there are two positive constants $C_i(d,q)$, $i=1,2$, such that for each $j\in [1,k]$ and any $f\in X_N^j$ we have
\be\label{G.3}
C_1(d,q)\|f\|_q^q \le \frac{1}{m} \sum_{\nu=1}^m |f(\xi^\nu)|^q \le C_2(d,q)\|f\|_q^q.
\ee
In the case $q=\infty$  for each $j\in [1,k]$ and any $f\in X_N^j$ we have
\be\label{G.4}
C_1(d)\|f\|_\infty \le \max_{1\le\nu\le m} |f(\xi^\nu)| \le  \|f\|_\infty.
\ee

In \cite{VT160} we studied the universal discretization  for the collection of subspaces of trigonometric polynomials with frequencies from parallelepipeds (rectangles). For $\bs\in\N^d_0$
define
$$
R(\bs) := \{\bk \in \Z^d :   |k_j| < 2^{s_j}, \quad j=1,\dots,d\}.
$$
    Let $Q$ be a finite subset of $\Z^d$. We denote
$$
\Tr(Q):= \{f: f=\sum_{\bk\in Q}c_\bk e^{i(\bk,\bx)}\}.
$$
Consider the collection $\cC(n,d):= \{\Tr(R(\bs)), \|\bs\|_1=n\}$.

The following theorem was proved in \cite{VT160}.
\begin{Theorem}\label{GT.1} Let a set $T$ with cardinality $|T|= 2^r=:m$ have dispersion 
satisfying the bound disp$(T) < C(d)2^{-r}$ with some constant $C(d)$. Then there exists 
a constant $c(d)\in \N$ such that the set $2\pi T:=\{2\pi\bx: \bx\in T\}$ provides the universal discretization in $L_\infty$ for the collection $\cC(n,d)$ with $n=r-c(d)$.
\end{Theorem}
Theorem \ref{GT.1} is a conditional result. As we discussed in Section \ref{F} existence of sets with a property required in Theorem \ref{GT.1} is a non-trivial fact. In particular, the $(t,r,d)$-nets provide such existence. We now give a definition of the $(t,r,d)$-nets.
\begin{Definition}\label{GD.1} A $(t,r,d)$-net (in base $2$) is a set $T$ of $2^r$ points in 
$[0,1)^d$ such that each dyadic box $[(a_1-1)2^{-s_1},a_12^{-s_1})\times\cdots\times[(a_d-1)2^{-s_d},a_d2^{-s_d})$, $1\le a_j\le 2^{s_j}$, $j=1,\dots,d$, of volume $2^{t-r}$ contains exactly $2^t$ points of $T$.
\end{Definition}
We note that existence of $(t,r,d)$-nets is a very non-trivial problem.  A construction of such nets for all $d$ and $t\ge Cd$, where $C$ is a positive absolute constant, $r\ge t$  is given in \cite{NX}.

Theorem \ref{GT.1} in a combination with Theorems \ref{FT.1} and \ref{FT.2} guarantees that the appropriately chosen Fibonacci ($d=2$) and  Frolov (any $d\ge 2$) point sets provide 
universal discretization in $L_\infty$ for the collection $\cC(n,d)$. 

The following Theorem \ref{GT.1'} (see \cite{VT160}) can be seen as an inverse to Theorem~\ref{GT.1}.

 \begin{Theorem}\label{GT.1'}  Assume that $T\subset [0,1)^d$ is such that the set $2\pi T$ provides universal discretization in $L_\infty$ for the collection
 $\cC(n,d)$ with a constant $C_1(d)$ (see (\ref{G.2})). Then there exists a positive constant $C(d)$ with the following property disp$(T) \le C(d)2^{-n}$.
 \end{Theorem}

{\bf Arbitrary trigonometric polynomials.}  For $n\in \N$ denote $\Pi_n :=\Pi(\bN)\cap \Z^d$ with $\bN =(2^{n-1}-1,\dots,2^{n-1}-1)$, where $\Pi(\bN) := [-N_1,N_1]\times\cdots\times[-N_d,N_d]$.
 Then $|\Pi_n| = (2^n-1)^d <2^{dn}$. Let $v\in\N$ and $v\le |\Pi_n|$. Consider
 $$
 \cS(v,n):= \{Q\subset \Pi_n : |Q|=v\}.
 $$
 Then it is easy to see that
 \be\label{G.5}
 |\cS(v,n)| =\binom{|\Pi_n|}{v}<2^{dnv}.
 \ee

 We are interested in solving the following problem of universal discretization.
 For a given $\cS(v,n)$ and $q\in [1,\infty)$ find a condition on $m$ such that there exists a set
 $\xi = \{\xi^\nu\}_{\nu=1}^m$ with the property: for any $Q\in \cS(v,n)$ and each
 $f\in \Tr(Q)$ we have
 \be\label{G.6}
 C_1(q,d)\|f\|_q^q \le \frac{1}{m}\sum_{\nu=1}^m |f(\xi^\nu)|^q \le C_2(q,d)\|f\|^q_q.
 \ee

We present results from \cite{DPTT} for $q=2$ and $q=1$.

\begin{Theorem}\label{GT.2} There exist three positive constants $C_i(d)$, $i=1,2,3$,
such that for any $n,v\in\N$ and $v\le |\Pi_n|$ there is a set $\xi =\{\xi^\nu\}_{\nu=1}^m \subset \T^d$, with $m\le C_1(d)v^2n$, which provides universal discretization
in $L_2$ for the collection $\cS(v,n)$: for any $f\in \cup_{Q\in \cS(v,n)} \Tr(Q)$
$$
C_2(d)\|f\|_2^2 \le \frac{1}{m} \sum_{\nu=1}^m |f(\xi^\nu)|^2 \le C_3(d)\|f\|_2^2.
$$
\end{Theorem}

The classical Marcinkiewicz-type result for $\Tr(\Pi_n)$ provides a universal set $\xi$ with cardinality $m\le C(d)2^{dn}$. Thus, Theorem \ref{GT.2} gives a non-trivial result
for $v$ satisfying $v^2n\le C(d)2^{dn}$.

\begin{Theorem}\label{GT.3} There exist three positive constants $C_1(d)$, $C_2$, $C_3$,
such that for any $n,v\in\N$ and $v\le |\Pi_n|$ there is a set $\xi =\{\xi^\nu\}_{\nu=1}^m \subset \T^d$, with $m\le C_1(d)v^2n^{9/2}$, which provides universal discretization
in $L_1$ for the collection $\cS(v,n)$: for any $f\in \cup_{Q\in \cS(v,n)} \Tr(Q)$
$$
C_2\|f\|_1 \le \frac{1}{m} \sum_{\nu=1}^m |f(\xi^\nu)| \le C_3\|f\|_1.
$$
\end{Theorem}

The classical Marcinkiewicz-type result for $\Tr(\Pi_n)$ provides a universal set $\xi$ with cardinality $m\le C(d)2^{dn}$. Thus, Theorem \ref{GT.3} gives a non-trivial result
for $v$ satisfying $v^2n^{9/2}\le C(d)2^{dn}$.
 
\section{Generalizations}
\label{H}

As above for a function class $\bW$ we have a concept of error of the cubature formula $\La_m(\cdot,\xi)$
\be\label{H.3}
\La_m(\bW,\xi):= \sup_{f\in \bW} |\int_\Omega fd\mu -\La_m(f,\xi)|. 
\ee
  If the class $\bW=\{f(\bx,\by): \by \in Y\}$ is parametrized by a parameter $\by \in Y\subset \R^n$ with $Y$ being a bounded measurable set, then we can consider a natural {\it average case setting}. For $\bp =(p_1,\dots,p_n)$ define
\be\label{H.4}
\La_m(\bW,\xi,\bp):= \|\int_\Omega f(\cdot,\by)d\mu -\La_m(f(\cdot,\by),\xi)\|_\bp,
\ee
where the vector $L_\bp$ norm is taken with respect to the Lebesgue measure on $Y$. 
We write $\La_m(\bW,\xi,\infty):= \La_m(\bW,\xi)$. We are interested in dependence on $m$ of the quantities
$$
\kappa_m (\bW,\bp) := \inf_{ \lambda_1,\dots,\lambda_m; \xi^{1},\dots,
\xi^m}\Lambda_m(\bW,\xi,\bp)
$$
for different classes $\bW$.

 We now present a rather general setting of this problem.  Let $1\le q\le \infty$. We define a set $\mathcal K_q$ of kernels possessing the following properties. 
 Let $K(\bx,\by)$ be a measurable function on $\Omega^1\times\Omega^2$.
 We assume that for any $\bx\in\Omega^1$ we have $K(\bx,\cdot)\in L_q(\Omega^2)$; for any $\by\in \Omega^2$ the $K(\cdot,\by)$ is integrable over $\Omega^1$ and $\int_{\Omega^1} K(\bx,\cdot)d\bx \in L_q(\Omega^2)$. For $1\le p\le \infty$ and a kernel $K\in \mathcal K_{p'}$, $p':=p/(p-1)$, we define the class
\be\label{H.5}
\bW^K_p :=\{f:f=\int_{\Omega^2}K(\bx,\by)\varphi(\by)d\by,\quad\|\varphi\|_{L_p(\Omega^2)}\le 1\}.  
\ee
Then each $f\in \bW^K_p$ is integrable on $\Omega^1$ (by Fubini's theorem) and defined at each point of $\Omega^1$. We denote for convenience
$$
 J_K(\by):=\int_{\Omega^1}K(\bx,\by)d\bx.
$$

For a cubature formula $\Lambda_m(\cdot,\xi)$ we have
$$
\Lambda_m(\bW^K_p,\xi) = \sup_{\|\varphi\|_{L_p(\Omega^2)}\le 1} |\int_{\Omega^2}\bigl( J_K(\by)-\sum_{\mu=1}^m\lambda_\mu K(\xi^\mu,\by)\bigr)\varphi(\by)d\by|=
$$
\be\label{H.6}
=\| J_K(\cdot)-\sum_{\mu=1}^m\lambda_\mu K(\xi^\mu,\cdot)\|_{L_{p'}(\Omega^2)}.
\ee
Consider a problem of numerical integration of functions  $K(\bx,\by)$, $\by\in\Omega^2$, with respect to  $\bx$,  $K\in {\mathcal K}_q$, in other words a problem of numerical integration of functions from the function class $\bK:=\{K(\bx,\by):\by\in\Omega^2\}$:
$$
{ \int_{\Omega^1} K(\bx,\by)d\bx - \sum_{\mu=1}^m \lambda_\mu K(\xi^\mu,\by)}.
$$

\begin{Definition}\label{HD.1}  $(K,q)$-discrepancy of a set of 
knots  $\xi^1,\dots,\xi^m$ and a set of weights  $\lambda_1,\dots,\lambda_\mu$ (a cubature formula  $(\xi,\Lambda)$) is
$$
 D(\xi,\Lambda,K,q):=\Lambda_m(\bK,\xi,q)=\|\int_{\Omega^1} K(\bx,\by)d\bx - \sum_{\mu=1}^m \lambda_\mu K(\xi^\mu,\by)\|_{L_q(\Omega^2)}.
$$
\end{Definition} 
In a special case $\La_m(\cdot,\xi) = Q_m(\cdot,\xi)$ we write $D(\xi,Q,K,q)$.
The above definition of the $(K,q)$-discrepancy and relation (\ref{H.6})  imply right a way the following relation
\be\label{H.7}
D(\xi,\Lambda,K,p') = \Lambda_m(\bW^K_p,\xi).
\ee
Relation (\ref{H.7}) shows that numerical integration in the class $\bW^K_p$ and the $(K,q)$-discrepancy are 
tied by the duality principle.

Let us consider a special case, when $K(\bx,\by)=F(\bx-\by)$, $\Omega^1=\Omega^2 = [0,1)^d$ and we deal with $1$-periodic in each variable functions. 
Associate with a cubature formula $(\xi,\La)$ and the function $F$ the following function
$$
g_{\xi,\Lambda,F}( \bx) := \sum_{ \bk}\Lambda(\xi,\bk)\hat F
( \bk)e^{2\pi i( \bk, \bx)}- \hat F(\mathbf 0),
$$
where
$$
\Lambda(\xi,\bk):= \La_m(e^{2\pi i(\bk,\bx)},\xi).
$$
Then for the quantity $\Lambda_m ( \bW_{p}^F ,\xi)$
we have ($p' := p/(p-1)$)
$$
\Lambda_m ( \bW_{p}^F ,\xi)=
\sup_{f\in \bW_{p}^F}
\bigl| \Lambda_m(f,\xi) -\hat f(\mathbf 0)\bigr|
$$
$$
=\sup_{\|\varphi\|_p\le 1} \bigl| \Lambda_m \bigl(F( \bx)\ast
\varphi( \bx),\xi\bigr) -\hat F(\mathbf 0)\hat\varphi(\mathbf 0)\bigr| 
$$
\be\label{H.8}
=\sup_{\|\varphi\|_p\le 1}\bigl|\<g_{\xi,\Lambda,F}(-\by),
\overline{\varphi(\by)}\>\bigr|=\|g_{\xi,\Lambda,F}\|_{p'}.
\ee
Let us discuss a special case of function $F$, which is very important in numerical
integration (see, for instance, \cite{TBook}, \cite{VT89}, \cite{VTbookMA}, and \cite{DTU}).  Let for $r>0$ 
 \be\label{H.9}
F_{r,\alpha}(x):= 1+2\sum_{k=1}^\infty k^{-r}\cos (2\pi kx-\alpha \pi/2).
\ee
For $\bx=(x_1,\dots,x_d)$, $\alpha=(\alpha_1,\dots,\alpha_d)$ denote
$$
F_{r,\alpha}(\bx) := \prod_{j=1}^d F_{r,\alpha_j}(x_j)
$$
and
$$
\bW^r_{p,\alpha} :=\bW^{F_{r,\alpha}}_p= \{f:f(\bx)=(F_{r,\alpha}\ast \varphi)(\bx)
$$
$$
:= \int_{[0,1)^d} F_{r,\alpha}(\bx-\by)\varphi(\by)d\by,\quad \|\varphi\|_p \le 1\}.
$$ 
  In the case of integer $r$ the class $\bW^r_{p,\alpha}$ with $\alpha=(r,\dots,r)$ is very close to the class of functions $f$, 
satisfying $\|f^{(r,\dots,r)}\|_p \le 1$, where $f^{(r,\dots,r)}$ is the mixed derivative of $f$ of order $rd$. 

It is easy to see that
\be\label{H.10}
\|g_{\xi,\Lambda,F_{r,\alpha}}\|_{2} = \left(\sum_{\bk\neq \mathbf 0}|\Lambda(\xi,\bk)|^2\left(\prod_{j=1}^d (\max(|k_j|,1))^{-r}\right)^2 + |\La(\xi,\mathbf 0)-1|^2\right)^{1/2} .
\ee
The above quantity in the case $r=1$ was introduced in \cite{Zi} under the name 
{\it diaphony}. In case of generic $r$ it was called {\it generalized diaphony} and was studied in \cite{Lev7}. Relation (\ref{H.8}) shows that generalized diaphony is closely related to numerical integration of the class $\bW^r_{2,\alpha}$. Following this analogy,
we can call the quantity $\|g_{\xi,\Lambda,F_{r,\alpha}}\|_{q}$ the $(r,q)$-diaphony of the pair $(\xi,\La)$ (the cubature formula $(\xi,\La)$). Behavior of $\kappa_m(\bW^r_{p,\alpha})$ is well studied (see, for instance, \cite{VTbookMA} and \cite{DTU}). By (\ref{H.8}) results on 
$\kappa_m(\bW^r_{p,\alpha})$ provide estimates on
$$
\inf_{\xi,\Lambda}\|g_{\xi,\Lambda,F_{r,\alpha}}\|_{p'}.
$$

For completeness we cite some known results on the lower bounds for $\kappa_m(\bW^r_{p,\alpha})$. The reader can find these and other results with a historical discussion in \cite{VTbookMA}, Chapter 6 and in \cite{DTU}, Chapter 8. 

\begin{Theorem}\label{HT.1} The following lower estimate is valid for any
cubature formula $(\xi,\Lambda)$ with $m$ knots $(r > 1/p)$
$$
\Lambda_m( \bW_{p,\alpha}^r,\xi) \ge C(r,d,p)m^{-r}
(\log m)^{\frac{d-1}{2}},\qquad 1 \le p < \infty .
$$
\end{Theorem}
The rate of decay $m^{-r}(\log m)^{\frac{d-1}{2}}$ in the lower bound in Theorem \ref{HT.1} does not depend on $p$. Therefore, the larger the $p<\infty$ the stronger the lower bound. 
It turns out that in the case $p=1$ one can improve the corresponding lower bound under certain restrictions on the weights of the cubature formula. 
We obtained the lower estimates for the quantities
$$
     \kappa_m^B(\bW) := \inf_{\Lambda_m(\cdot,\xi) \in Q(B,m)} \Lambda_m(\bW,\xi).
$$
 We proved the following relation.      
 \begin{Theorem}\label{HT.2} Let $r>1$. Then
$$
     \kappa_m^B(\bW_{1,0}^r) \geq   C(r,B,d)   m^{-r}(\log   m)^{d-1},   \qquad
C(r,B,d)>0.
$$
\end{Theorem}

The case $p=\infty$ is excluded in Theorem \ref{HT.1}. There is no nontrivial general lower estimates in this case. We give one conditional result in this direction.

\begin{Theorem}\label{HT.3} Let the cubature formula $(\xi,\Lambda)$
be such that the inequality
$$
\Lambda_m ( \bW_{p,\alpha}^r,\xi)\le C_1(p,r,d) m^{-r} (\log m)^{(d-1)/2},
\qquad r > 1/p,
$$
holds for some $1 < p < \infty$.
Then there exists a constant $C_2(p,r,d)>0$ such that 
$$
 \Lambda_m ( \bW_{\infty,\alpha}^r ,\xi)\ge C_2(p,r,d) m^{-r}
(\log m)^{(d-1)/2} .
$$
\end{Theorem}

There are two big open problems in this area. We formulate them as conjectures.
\begin{Conjecture}\label{HC.1} For any $d\ge 2$ and any $r\ge 1$ we have
$$
\kappa_m(\bW^r_{1,\alpha}) \ge C(r,d)m^{-r}(\log m)^{d-1}.
$$
\end{Conjecture}
\begin{Conjecture}\label{HC.2} For any $d\ge 2$ and any $r> 0$ we have
$$
\kappa_m(\bW^r_{\infty,\alpha}) \ge C(r,d)m^{-r}(\log m)^{(d-1)/2}.
$$
\end{Conjecture}
We note that by   Theorem \ref{CT.1} and (\ref{C.2}) Conjecture \ref{HC.1} implies Conjecture \ref{BC.1} and Conjecture \ref{HC.2} implies for any cubature formula $(\xi,\Lambda)$
\be\label{H.11}
D^r_1(\xi,\Lambda) \ge C(r,d)m^{-r}(\log m)^{(d-1)/2}.  
\ee
 
 \begin{Remark}\label{HR.1} In the case $d=2$, $r=1$, and $\alpha=(1,1)$ Conjecture \ref{HC.1} holds.
 \end{Remark}
 
 Remark \ref{HR.1} follows from an analog of the Schmidt's bound (\ref{B.2}) and Proposition \ref{CP.1}. We discuss this in more detail. D. Bilyk and I observed that a slight modification of 
 the proof of (\ref{B.2}) from \cite{Bi} gives the following lower bound. For any cubature formula $(\xi,\Lambda)$ we have
 \be\label{H.12}
 \Lambda_m(\chi^2,\xi) \ge C_1m^{-1}\log m.
 \ee
 Therefore, by Proposition \ref{CP.1} for any cubature formula $(\xi,\Lambda)$, satisfying an extra condition $\sum_j \lambda_j =1$, we have for $\bW^1_1 := \bW^1_{1,(1,1)}$
 \be\label{H.13}
 \Lambda_m(\bW^1_1,\xi) \ge C_2m^{-1}\log m.
 \ee
 Further, it is well known (see \cite{VTbookMA}, p.269) and easy to check, that for a function class $\bW$ of periodic functions, satisfying the condition: $1\in \bW$ and for $f\in \bW$ we have $\frac{1}{2}(f-\hat f(\mathbf 0))\in \bW$, the inequality holds
\be\label{H.14}
 \inf_\Lambda\Lambda_m(\bW,\xi) \ge \frac{1}{4}  \inf_{\Lambda:\sum_j \lambda_j =1}\Lambda_m(\bW,\xi).
 \ee
 Clearly, $\bW^1_1$ satisfies the above condition on a class $\bW$. Combining (\ref{H.12})--(\ref{H.14}) we obtain for $d=2$
$$
\kappa_m(\bW^1_{1}) \ge Cm^{-1}\log m.
$$

\section{Numerical integration without smoothness assumptions}
\label{J}

In the previous sections we discussed numerical integration for 
classes of functions under certain conditions on smoothness. Parameter $r$ controlled 
the smoothness. The above results show that the numerical integration characteristics 
decay with the rate $m^{-r}(\log m)^{c(d)}$, which substantially depends on smoothness $r$. The larger the smoothness -- the faster the error decay. In this section 
we discuss the case, when we do not impose any of the smoothness assumptions. 
Surprisingly, even in such a situation we can guarantee some rate of decay. Results discussed in this section apply in a very general setting. We present here results from \cite{VT164}.
The following result is proved 
in \cite{VT149} (see also \cite{VT89} for previous results). For the theory of greedy algorithms we refer the reader to \cite{Tbook}.
Consider a dictionary
$$
\Di := \{K(\bx,\cdot), \bx\in \Omega^1\}
$$
and define a Banach space $X(K,q)$ as the $L_{q}(\Omega^2)$-closure of span of $\Di$. 

\begin{Theorem}\label{JT.1} Let { $\bW^K_p$} be a class of functions defined above in Section \ref{H}. Assume that { $K\in \mathcal K_{p'}$} satisfies the condition
$$
 { \|K(\bx,\cdot)\|_{L_{p'}(\Omega^2)} \le 1, \quad \bx\in \Omega^1,\quad |\Omega^1|=1},
$$
and { $J_K\in X(K,p')$}. Then for any { $m$} there exists (provided by an appropriate greedy algorithm)  a cubature formula { $Q_m(\cdot,\xi)$} such that
$$
 Q_m(\bW^K_p,\xi)\le C(p-1)^{-1/2} m^{-1/2}, \quad 1< p\le 2 .
$$
\end{Theorem}
As a direct corollary of Theorem \ref{JT.1} and relation (\ref{H.7}) we obtain the following result about the $(K,q)-discrepancy$. 
  
\begin{Theorem}\label{JT.2}   Assume that { $K\in \mathcal K_{q}$} satisfies the condition
$$
 { \|K(\bx,\cdot)\|_{L_{q}(\Omega^2)} \le 1, \quad \bx\in \Omega^1,\quad |\Omega^1|=1},
$$
and  $J_K\in X(K,q)$. Then for any  $m$ there exists (provided by an appropriate greedy algorithm)  a cubature formula  $Q_m(\cdot,\xi)$  such that
$$
 D(\xi,Q,K,q)\le Cq^{1/2} m^{-1/2}, \quad 2\le q <\infty.
$$
\end{Theorem}

\begin{Remark}\label{JR.1} In Theorems \ref{JT.1} and \ref{JT.2} we impose the restriction $1<p\le 2$ or the dual one $2\le q <\infty$. The proof of Theorems \ref{JT.1} and \ref{JT.2} from \cite{VT149} also works in the case $2<p<\infty$ or $1<q<2$ and gives
$$
 Q_m(\bW^K_p,\xi)\le C  m^{-1/p}, \quad 2< p<\infty ,
$$
$$
 D(\xi,Q,K,q)\le C m^{\frac{1}{q}-1}, \quad 1< q <2.
$$
\end{Remark}

Let us discuss a special case $K(\bx,\by)=F(\bx-\by)$, $\Omega^1=\Omega^2 = [0,1)^d$ and $1$-periodic in each variable functions. Then we associate with a cubature formula $(\xi,\La)$ and the function $F$ the function $g_{\xi,\La,F}(\bx)$.
The following Proposition is proved in \cite{VT89}.
\begin{Proposition}\label{JP.1} Let $1<p<\infty$ and $\|F\|_{p'} \le 1$. Then the kernel
$K(\bx,\by)=F(\bx-\by)$ satisfies the assumptions of Theorem \ref{JT.1}.
\end{Proposition}
Proposition \ref{JP.1}, Theorem \ref{JT.1}, Remark \ref{JR.1}, and relation (\ref{H.8}) imply
\begin{Theorem}\label{JT.3} Let $1<p<\infty$ and $\|F\|_{p} \le 1$. Then there exists 
a set $\xi$ of $m$ points such that
$$
\|g_{\xi,Q,F}(\bx)\|_p \le Cp^{1/2} m^{-1/2}, \quad 2\le p<\infty, 
$$
$$
\|g_{\xi,Q,F}(\bx)\|_p \le C m^{\frac{1}{p}-1}, \quad 1< p <2.
$$
\end{Theorem}

Here is a corollary of Theorem \ref{JT.2} and Proposition \ref{JP.1}. 
Let $E\subset [0,1)^d$ be a measurable set. Consider $F(\bx):=\tilde \chi_E(\bx)$. 
\begin{Theorem}\label{JT.4} For any $p\in [2,\infty)$ there exists a set of $m$ points $\xi$ such that  
$$
Q_m(\{\tilde \chi_E(\bx-\bz),\bz\in [0,1)^d\},\xi,p) \le Cp^{1/2} m^{-1/2}.
$$
\end{Theorem}
We note that there are interesting results on the behavior of $Q_m(\{\chi_E(\bx-\bz),\bz\in [0,1)^d\},\xi,\infty)$ under assumption that $E$ is a convex set (see \cite{BC}). 
Theorem \ref{JT.4} shows that for $p<\infty$ we do not need any assumptions on the geometry of $E$ in order to get the upper bound $\ll m^{-1/2}$ for the discrepancy. 

The proof of the above Theorems \ref{JT.1}--\ref{JT.4} is constructive (see \cite{VT149}), it is based on the greedy algorithms. We formulate the related result from the theory of greedy approximation.   We remind some notations from the theory of greedy approximation in Banach spaces. The reader can find a systematic presentation of this theory in \cite{Tbook}, Chapter 6. 
Let $X$ be a Banach space with norm $\|\cdot\|$. We say that a set of elements (functions) ${\mathcal D}$ from $X$ is a dictionary   if each $g\in {\mathcal D}$ has norm less than or equal to one ($\|g\|\le 1$)
 and the closure of $\sp {\mathcal D}$ coincides with $X$.  
 
For an element $f\in X$ we denote by $F_f$ a norming (peak) functional for $f$: 
$$
\|F_f\| =1,\qquad F_f(f) =\|f\|.
$$
The existence of such a functional is guaranteed by the Hahn-Banach theorem.

We proceed to  the Incremental Greedy Algorithm (see \cite{T12} and \cite{Tbook}, Chapter 6).       Let $\ep=\{\ep_n\}_{n=1}^\infty $, $\ep_n> 0$, $n=1,2,\dots$ . For a Banach space $X$ and a dictionary $\Di$ define the following algorithm IA($\ep$) $:=$ IA($\ep,X,\Di$).

 {\bf Incremental Algorithm with schedule $\ep$ (IA($\ep,X,\Di$)).} 
  Denote $f_0^{i,\ep}:= f$ and $G_0^{i,\ep} :=0$. Then, for each $m\ge 1$ we have the following inductive definition.

(1) $\ff_m^{i,\ep} \in \Di$ is any element satisfying
$$
F_{f_{m-1}^{i,\ep}}(\ff_m^{i,\ep}-f) \ge -\ep_m.
$$

(2) Define
$$
G_m^{i,\ep}:= (1-1/m)G_{m-1}^{i,\ep} +\ff_m^{i,\ep}/m.
$$

(3) Let
$$
f_m^{i,\ep} := f- G_m^{i,\ep}.
$$

We consider here approximation in uniformly smooth Banach spaces. For a Banach space $X$ we define the modulus of smoothness
$$
\rho(u) := \sup_{\|x\|=\|y\|=1}\left(\frac{1}{2}(\|x+uy\|+\|x-uy\|)-1\right).
$$
It is well known that in the case $X=L_p$, 
$1\le p < \infty$ we have
\begin{equation}\label{J.1}
\rho(u) \le \begin{cases} u^p/p & \text{if}\quad 1\le p\le 2 ,\\
(p-1)u^2/2 & \text{if}\quad 2\le p<\infty. \end{cases}     
\end{equation}
 
 Denote by $A_1({\mathcal D}):=A_1(\Di,X)$ the closure in $X$ of the convex hull of ${\mathcal D}$. Proof of Theorem \ref{JT.1} and Remark \ref{JR.1} is based on the following theorem proved in \cite{T12} (see also \cite{Tbook}, Chapter 6).

\begin{Theorem}\label{JT.5} Let $X$ be a Banach space with  modulus of smoothness $\rho(u)\le \gamma u^q$, $1<q\le 2$. Set
$$
\ep_n := \bt\gamma ^{1/q}n^{-1/p},\qquad p:=\frac{q}{q-1},\quad n=1,2,\dots .
$$
Then, for every $f\in A_1({\mathcal D})$ we have
$$
\|f_m^{i,\ep}\| \le C(\bt) \gamma^{1/q}m^{-1/p},\qquad m=1,2\dots.
$$
\end{Theorem}

{\bf Acknowledgment.}   The work was supported by the Russian Federation Government Grant No. 14.W03.31.0031.

\end{document}